\documentclass{ws-ijgmmp}
\expandafter\let\csname equation*\endcsname\relax
\expandafter\let\csname endequation*\endcsname\relax 
\usepackage{amsmath,amsfonts,amssymb}
\usepackage{calc,graphicx,color}
\usepackage{hyperref}

\newcommand{\figtobox}[3][{}]{%
	\expandafter\newsavebox\csname #2box\endcsname%
	\expandafter\savebox\csname #2box\endcsname{#1\input{#3.latex}}%
	\expandafter\def\csname #2\endcsname%
		{\setlength{\figoff}{1.0ex-.5\expandafter\ht\csname #2box\endcsname}%
			\raisebox{\figoff}{\expandafter\usebox\csname #2box\endcsname}}%
}
\newlength{\figoff}

\figtobox{jbasis}{jbasis}
\figtobox{lbasis}{lbasis}
\figtobox{lnorm}{lnorm}
\figtobox{jnorm}{jnorm}
\figtobox{Llmat}{Llmat}
\figtobox{Ljmat}{Ljmat}
\figtobox{Lvert}{Lvert}
\figtobox{bblob}{bblob}
\figtobox{bblobo}{bblobo}
\figtobox{bblobv}{bblobv}
\figtobox[]{blob}{blob}
\figtobox[]{blobo}{blobo}
\figtobox[]{blobv}{blobv}
\figtobox[]{thabc}{theta}
\figtobox[]{thbcj}{theta}
\figtobox[]{thadjp}{theta}
\figtobox[]{thadj}{theta}
\figtobox[]{thabl}{theta}
\figtobox[]{thcdl}{theta}
\figtobox[]{thjjp}{theta}
\figtobox[]{thjjpp}{theta}
\figtobox[]{thjjptwo}{theta}
\figtobox[]{loopa}{loop}
\figtobox[]{loopl}{loop}
\figtobox[]{loopj}{loop}
\figtobox[]{loopjp}{loop}
\figtobox[]{tetjl}{tet}
\figtobox[]{tetll}{tet}
\figtobox[]{tetaadjj}{tet}
\figtobox[]{tetbbcjj}{tet}
\figtobox[]{tetadjj}{tet}

\figtobox{blobss}{blob3}
\figtobox{blobsss}{blob4}

\def\braket#1#2{\langle #1 | #2 \rangle}
\def\bra#1{\langle#1|}
\def\ket#1{|#1\rangle}
\def\Tet{\operatorname{Tet}}
\newcommand{\C}{\mathbb{C}}
\newcommand{\sixjkl}[6]{\begin{Bmatrix}{#1}&{#2}&{#3}\\{#4}&{#5}&{#6}\end{Bmatrix}_{KL}}
\newcommand{\sixjrw}[6]{\begin{Bmatrix}{#1}&{#2}&{#3}\\{#4}&{#5}&{#6}\end{Bmatrix}_{RW}}
\newcommand{\ul}{\underline{l}}
\newcommand{\uj}{\underline{j}}
\newcommand{\ol}{\overline{l}}
\newcommand{\oj}{\overline{j}}
\newcommand{\su}{\mathrm{su}}

\begin{document}
%
%
\catchline{}{}{}{}{}
%
\title{RECURRENCE RELATION FOR THE $6j$-SYMBOL OF $\su_q(2)$ AS A SYMMETRIC EIGENVALUE PROBLEM}
\author{\footnotesize IGOR KHAVKINE}
\address{Department of Mathematics,
University of Trento,
I--38123 Povo (TN) Italy\\
TIFPA-INFN, Trento\\
\email{igor.khavkine@unitn.it}}

\maketitle

\begin{history}
\received{(Day Month Year)}
\revised{(Day Month Year)}
\accepted{(Day Month Year)}
\comby{(xxxxxxxxx)}
\end{history}

\begin{abstract}
A well known recurrence relation for the $6j$-symbol of the quantum
group $\su_q(2)$ is realized as a tridiagonal, symmetric eigenvalue
problem. This formulation can be used to implement an efficient
numerical evaluation algorithm, taking advantage of existing specialized
numerical packages. For convenience, all formulas relevant for such an
implementation are collected in the appendix. This realization is a
byproduct of an alternative proof of the recurrence relation, which
generalizes a classical ($q=1$) result of Schulten and Gordon and uses
the diagrammatic spin network formalism of Temperley-Lieb recoupling
theory to simplify intermediate calculations.
\end{abstract}

\keywords{Quantum group; $6j$-symbol; recurrence relation; eigenvalue problem.}

\ccode{2010 Mathematics Subject Classification: 20G05, 20G42, 65Q30, 65F15}

\section{Introduction}
Quantum groups first appeared in the study of quantum integrable
systems. Since then, they have proven useful in many applications,
including among others conformal field theory, statistical mechanics,
representation theory and the theory of hypergeometric functions, along
with exhibiting a rich internal structure. Quantum groups, and
$\su_q(2)$ in particular, have appeared seminally in the mathematical
physics literature in connection with topological quantum field
theory~\cite{witten}, notably in the Turaev-Viro~\cite{TV}
``regularized'' combinatorial model of 3d quantum gravity, inspired by
the $\su(2)$-based model of Ponzane and Regge~\cite{PR}. More recently,
the quantum group $\su_q(2)$ was used to construct
``loop''~\cite{smolin2,BMS,smolin} and ``spin foam''~\cite{KC,FM,DMBS}
models of quantum gravity with a positive cosmological constant. Many of
these applications crucially involve the $\su(2)$ and $\su_q(2)$
$6j$-symbols and identities associated with them.

Not surprisingly, the $6j$-symbol (or Racah-Wigner coefficients and
henceforth the $6j$ or the $q$-$6j$ for the classical and $q$-deformed
versions respectively) has been studied extensively. It first appeared
in work on $q$-hypergeometric functions~\cite{AW-rec}. Later, it was
found to play a central role in the representation theory of
$\su_q(2)$~\cite{KR-q6j}. It is known to satisfy some recurrence
relations~\cite{AW-rec,KK-rec}, including a particular linear, three
term, single argument one~\cite{KR-q6j,MT-ptps,MT-prl}. This recurrence
has been used to analyze the asymptotics of the
classical~\cite{PR,SG-rec} and quantum~\cite{WT} $6j$.

One of the main goals of this paper is to give a new, simpler derivation
of this important recurrence relation for the $q$-$6j$, based on
graphical manipulations in the so-called \emph{spin network} formalism.
The second main goal of this paper, besides shining a spotlight on this
recurrence relation, is to show that the new derivation actually
realizes the recurrence as a tridiagonal, symmetric eigenvalue problem.
The latter property is special, in the sense of not being shared by most
recurrence relations, and appears to have been missed in the existing
literature.

Even though an explicit formula for the $q$-$6j$ is known, it is
important to spot-light the existence the above recurrence relation and
its properties for the following reasons. There exist segments of the
literature~\cite{KC,RRS} that require the numerical evaluation of
$q$-$6j$ symbols for large numbers of arguments, yet are not as well
aware of the existence of this recurrence relation as the more abstract
literature on quantum groups cited earlier. The explicit formula for the
$q$-$6j$ involves a number of arithmetic operations that is linear in
its arguments (see Appendix). In applications where a large number of
$6j$-symbols is needed at once, e.g.,\ for all values of one argument
with others fixed as is the case in~\cite{KC}, the total number of
operations becomes quadratic in the arguments. Using the recurrence
relation instead can greatly increase the efficiency of the calculation
by reducing the total operation count to be linear in the arguments.
Moreover, a symmetric eigenvalue formulation has the significant
advantage of the possibility to make use of readily available, robust
linear algebra packages, such as LAPACK~\cite{Lap}, which automatically
take care of the important issues of numerical accuracy and stability.
When $q=1$ or when $q$ is a primitive root of unity, the inner product,
with respect to which the problem is symmetric, even becomes either
positive- or negative-definite and standard, specialized numerical
methods can be exploited to increase the efficiency of the calculation
even further. The final goal of this paper is to concisely collect all
the relevant information needed to readily implement such an efficient
$q$-$6j$ numerical evaluation algorithm without intimate familiarity
with the literature on quantum groups or $q$-hypergeometric functions.
It is also worth mentioning that a non-negligible advantage of the
simplicity of our derivation is that its steps can be easily followed by
a reader with minimal background in quantum groups, which is an ability
that is crucial for checking the correctness of any computer
implementation of the recurrence.

Sections~\ref{sec:spinnet} and~\ref{sec:ident}, which can be skipped by
those familiar with the mathematical literature on $\su_q(2)$ recoupling
theory, introduce the basic notions of the spin network
formalism~\cite{Pen,KL,CFS}, define the Kauffman-Lins convention for the
$6j$-symbol and summarize basic diagrammatic identities needed for
Section~\ref{sec:rec}, where the recurrence relation is realized as an
eigenvalue problem. This is accomplished as a byproduct of the
new proof of the recurrence itself that generalizes the
classical argument from the Appendix of~\cite{SG-rec}, where the
eigenvalue formulation was explicitly used and exploited (see
also~\cite{AABDR}). The diagrammatic spin network formalism makes all
intermediate calculations easy to check and reproduce. The Appendix
conveniently summarizes all formulas needed for a direct computer
implementation of the recurrence-based evaluation of the $q$-$6j$,
including the connection between the Kauffman-Lins and Racah-Wigner
notational convention, which is traditionally used in physics.

\section{Spin networks}\label{sec:spinnet}
In a variety of physical and mathematical applications, one often encounters
tensor contraction expressions, which could be expressed in one of the
following forms:
\begin{equation}\label{tens-con}
	T^{(i)_a(m)_e}_{(k)_c}
		= A^{(i)_a (j)_b}_{(k)_c (l)_d} B_{(j)_b}^{(l)_d (m)_e} ,
	\quad\qquad
	\blobss ,
	\quad\qquad
	\blobsss .
\end{equation}
In the first form, $T$, $A$ and $B$ are invariant tensors, with each
index, say $(i)_a$, transforming under a representation, in this case
labeled $a$, of a group or an algebra. The application at hand usually
calls for evaluating $T$, or at least simplifying it. An extensive
literature on this subject exists for the classical group $SU(2)$ or its
Lie algebra, a subject known as \emph{angular momentum
recoupling}~\cite{VMKh,YLV}. It is well known that such tensor
contractions can be very efficiently expressed, manipulated, and
simplified using diagrams known as \emph{spin networks}~\cite{Pen},
which are illusrated in the second form above. Tensors correspond to
vertices, ore more generally blobs, while indices correspond to
edges, with an internal edge corresponding to a contraction and the
direction distinguising upper from lower indices. Extensions of these
techniques~\cite{KL,CFS} are also known for the \emph{quantum} (or
\emph{$q$-deformed}, since they depend on an arbitrary complex number
$q\ne0$) analogs, the quantum group $\su_q(2)$ or $U_q(su(2))$. Once
certain normalization, ordering and index raising-lowering conventions
have been made, it is no longer required to keep track of the positions
(directions) and labels of individual indices (edges), so that only the
representation label needs to be kept. Such a convention is a great
advantange especially in the $q$-deformed case, where otherwise many
non-trivial multiplicative factors need to be explicitly displayed.
These more economical spin networks are illustrated in the third form
above. The translation between the diagrammatic notation of the above
second and third forms is precisely given in~\cite{KL,CFS}, though see
also the Appendix~A of~\cite{CCK} for a concise presentation. The basics
of this economical diagrammatic formalism, as needed for the derivation
of the recurrence relation, are given in this and next sections. All
relevant formulas, including explicit spin network evaluations in terms
of \emph{quantum integers} are listed in the Appendix.

Thus, single spin networks are edge-labeled graphs%
	\footnote{Spin networks are actually ribbon graphs, but since all
	diagrams in this paper are planar, the ribbon structure can be added
	through blackboard framing.}, %
where each vertex has valence either $1$ (free index) or $3$
(Clebsch-Gordan intertwiner). General spin networks are formal linear
combinations of single spin networks. Edges attached to univalent
vertices are called \emph{free}. Spin networks without free edges are
called \emph{closed}. Conventionally, the labels are either half-integers
(\emph{spins}) of integers (\emph{twice-spins}), which correspond
to irreducible representations of $\su_q(2)$. Reference~\cite{KL}
labels all spin networks with twice-spins. Unless otherwise indicated,
all conventions in this paper follow~\cite{KL}. Two spin networks may
be equal even if not represented with identical labeled graphs. A
complete description of these identities are given in~\cite{KL}
and~\cite{CFS}; their study constitutes \emph{spin network recoupling}
and is what allows us to equate%
	\footnote{A down-to-earth guide to this correspondence, for the
	classical $q=1$ case, can be found in Appendix~A
	of~\cite{CCK}. Complete details with proofs can be found
	in~\cite{CFS}.} %
spin networks with $\su_q(2)$-invariant tensors and their contractions.

In this correspondence, each index of a tensor, transforming under an
irreducible representation, corresponds to a spin network edge, labeled
by the same representation (free indices correspond to free edges). In
particular, a closed spin network corresponds to a complex number. Spin
networks form a graded algebra over $\C$ (as do tensors). The grading is
given by the number of free edges (free indices) and the product is
diagrammatic juxtaposition (tensor product).

\section{Diagrammatic identities}\label{sec:ident}
The spin networks with $n$ free edges  with fixed labels
(\emph{$n$-valent} spin networks) form a linear space with a natural
bilinear form (or inner product). Suppose that the free edges are
ordered in some canonical way, then, given two spin networks, we can
reflect one of them in a mirror and connect the free edges in order. As
an illustration, consider the inner products constructed in
Equations~\eqref{norm-j} and~\eqref{norm-l} out of the spin networks
shown in Equation~\eqref{lb-jb}, where we have taken the reflection
axis to be vertical. The value of the resulting closed spin network
defines the bilinear form, which is symmetric and
non-degenerate~\cite{KL}.  We use the bra-ket notation for this inner
product $\braket{s'}{s}$, where $s$ and $s'$ are two spin networks. We
also let $\ket{s}$ stand for $s$ and $\bra{s'}$ for the reflection of
$s'$. The existence of an inner product allows the following identities,
whose proofs can be found in~\cite{KL}. For each identity, the
corresponding well known fact of $SU(2)$ representation theory is given.

The space of $2$-valent spin networks, with ends labeled $a$ and $b$, is
$1$-dimensional if $a=b$, and $0$-dimensional otherwise.  For
non-trivial dimension, the single edge gives a complete basis and
therefore the \emph{bubble identity}:
\begin{equation}\label{bub-id}
	\bblob ~ = \delta_{ab} \frac{\bblobo}{\loopa} ~~ \bblobv ~ .
\end{equation}
This identity the diagrammatic analog of Schur's lemma for intertwiners
between irreducible representations.

The space of $3$-valent spin networks, with ends labeled $a$, $b$ and
$c$, is also $1$-dimensional if the triangle
inequalities~\eqref{tri-ineq} and parity constraints~\eqref{parity} are
satisfied, and $0$-dimensional otherwise, if $q$ is generic. When $q$ is
a primitive root of unity, the dimension also vanishes whenever the
further $r$-boundedness constraint~\eqref{r-bound} is violated. In the
case of nontrivial dimension, the canonical trivalent vertex gives a
complete basis and therefore the \emph{vertex collapse identity}:
\begin{align}\label{theta-id}
	\blobo ~ &= \frac{\blob}{\theta(a,b,c)} ~~ \blobv ~~, \\
	\theta(a,b,c) &= \thabc .
\end{align}
The normalization of the vertex, the value of the $\theta$-network, is
evaluated in Equation~\eqref{theta-def}. This identity is the diagrammatic
analog of the uniqueness (up to normalization) of the Clebsch-Gordan
intertwiner.

Now, consider the space of $4$-valent networks with free edges labeled
$a$, $b$, $c$ and $d$. There are two natural bases, the vertical
$\bra{l}$ and the horizontal $\ket{\bar{j}}$: 
\begin{equation}\label{lb-jb}
	\lbasis \quad,\quad \jbasis ~~.
\end{equation}
The admissible ranges for $j$ and $l$, the dimension $n$ of this space,
and the conditions on $(a,b,c,d)$ under which $n>0$ are given by
Equations~\eqref{jmin} through~\eqref{abcdrbd}. The transition matrix
between the two bases is given by the so-called $\Tet$-network:
\begin{equation}
	\Tet(a,b,c,d;j,l) = \tetjl = \braket{\bar{j}}{l}.
\end{equation}
The coefficients expressing the vertical basis in terms of the
horizontal one define the \emph{$6j$-symbol}, which can be expressed in
terms of the $\Tet$-network:
\begin{align}\label{KL-6j}
	\ket{l} &= \sum_{j} \sixjkl{a}{b}{j}{c}{d}{l} \ket{\bar{j}}, \\
	\sixjkl{a}{b}{j}{c}{d}{l} &= \frac{\tetjl}{\thadj\thbcj}~\loopj~.
\end{align}
Note the subscript $KL$ for Kauffman-Lins, since this $6j$-symbol is
defined with respect to the conventions of~\cite{KL}. The
relation to the classical Racah-Wigner $6j$-symbol used in the physics
literature is given explicitly in Equation~\ref{RW-KL}.

\section{Recurrence relation for the $\Tet$-network}\label{sec:rec}
The identities given in the previous section allow an new,
elementary derivation of the three-term recurrence relation for the
$\Tet$-network, distinct from the standard one. The standard derivation
is given in~\cite{MT-ptps} and another is possible using the general
theory of recurrences for $q$-hypergeometric functions~\cite{KK-rec},
but neither directly yields the symmetric eigenvalue problem form.

It is easy to check, using the bubble identity, that both the vertical
and horizontal bases are orthogonal and that they are normalized as
\begin{align}
\label{norm-j}
	\braket{\bar{j}}{\bar{j}} &= \jnorm~ = \frac{\thbcj ~ \thadj}{\loopj}, \\
	\braket l l               &= \lnorm~ = \frac{\thabl ~ \thcdl}{\loopl}.
\label{norm-l}
\end{align}
Curiously, when these normalizations
are fully expanded using formulas from the Appendix,
they take
the form $(-)^\sigma P/Q$, where $P$ and $Q$ are products of positive
quantum integers. In both cases, $\sigma=(a+b+c+d)/2$, is an integer independent of
$j$ or $l$. When $q=1$ or when $q$ is a primitive root of unity, positive
quantum integers are positive real numbers. Hence the above inner
product is real and either positive- or negative-definite. On the other
hand, for arbitrary complex $q$, the normalizations~\eqref{norm-j}
and~\eqref{norm-l} can be essentially arbitrary complex numbers.

If we can find a linear operator $L$ that is diagonal in one basis, but
not in the other, then we can obtain $\braket{\bar{j}}{l}$ as matrix
elements of the diagonalizing transformation. Furthermore, if the
non-diagonal form of $L$ is tridiagonal, then the linear equations
defining $\braket{\bar{j}}{l}$ reduce to a three-term recurrence
relation.

We can construct such an operator by generalizing the argument for the
classical case, found in the Appendix of~\cite{SG-rec}.  For
brevity of notation, we introduce a special modified version of the
trivalent vertex:
\begin{equation}
	\Lvert .
\end{equation}
The unlabeled edge implicitly carries twice-spin $2$ and the bold dot
indicates the multiplicative factor of $[a]$.  Using it, we can define a
symmetric operator $L$. Its diagrammatic representation and its matrix
elements are given below.

The operator $L$ is diagonal in the $\ket{l}$-basis
and its matrix elements $L_{ll'} = \bra{l} L \ket{l'}$ are
\begin{align}
	L_{ll'} &= \Llmat~ = \frac{[a][b]}{\loopl}~\thcdl~\tetll~\delta_{ll'} \\
		&= \lambda(a,b,l) \braket{l}{l'} ,
\end{align}
with
\begin{equation}
	\lambda(a,b,l) =
		\frac{\left[\frac{a-b+l}{2}\right]\left[\frac{-a+b+l}{2}\right]
			-\left[\frac{a+b-l}{2}\right]\left[\frac{a+b+l}{2}+2\right]}{[2]} ,
\end{equation}
where we have evaluated $\Tet(a,a,b,b;l,2)$ as
\begin{equation} 
	\label{tet-diag}
	\tetll = \frac{\thabl}{[a][b]} \lambda(a,b,l) .
\end{equation}
This result may be obtained directly from Equation~\eqref{tet-formula}, where
the sum reduces to two terms, or from more fundamental
considerations~\cite{MV}.
In the limit, $q\to 1$, the eigenvalues simplify to
$\lambda(a,b,l)=\frac{1}{4}[l(l+2)-a(a+2)-b(b+2)]$, which shows that the
operator $L$ is closely related to the ``square of angular momentum'' in
quantum mechanics, which was used to obtain the classical version of
this recurrence relation~\cite{SG-rec}.

On the other hand, in the $\ket{\bar{j}}$ basis, the operator
$L$ is not diagonal and the matrix elements
$\bar{L}_{jj'}=\bra{\bar{j}}L\ket{\bar{j}'}$, making use of the vertex
collapse identity, are
\begin{equation}
	\bar{L}_{jj'} = \Ljmat~ = \frac{[a][b]}{\thjjpp}~\tetaadjj~\tetbbcjj,
\end{equation}
with the special case $\bar{L}_{00}=0$. Fortunately, though
$\bar{L}_{jj'}$ is not diagonal, it is tridiagonal. This property is a
consequence of the conditions enforced at the central vertex in both
$\Tet$-networks above: the triangle inequality, $|j-j'|\le 2$, and the
parity constraint, which forces admissible values of $j$ to change by
$2$. If these conditions are violated, the matrix element
$\bar{L}_{jj'}$ vanishes.

The diagonal elements $\bar{L}_{jj}$ can be evaluated
using~\eqref{tet-diag}.
For the off-diagonal elements $\bar{L}_{j+2,j}=\bar{L}_{j,{j+2}}$
we also need
\begin{equation} 
	\label{tet-offdiag}
	\tetadjj = \frac{1}{[a]}\left[\frac{a+d-j}{2}\right]~\thadjp,
\end{equation}
which can be obtained in the same way as~\eqref{tet-diag}. Finally,
we need the identities
\begin{equation} 
	\thjjp = -\frac{[j+2]}{[2][j]}~\loopj
	\quad \text{and} \quad
	\thjjptwo = \loopjp \phantom{j+{}}.
\end{equation}
The $\ket{j}$-basis matrix elements can now be expressed as (again,
recall the special case $\bar{L}_{00}=0$)
\begin{align}
	\bar{L}_{jj} &= -\braket{j}{j} \frac{[2]\lambda(a,j,d)\lambda(b,j,c)}{[j][j+2]}, \\
	\bar{L}_{j,j+2} &= \braket{j{+}2}{j{+}2} \left[\frac{a+d-j}{2}\right]
			\left[\frac{b+c-j}{2}\right].
\end{align}

The transition matrix elements $\braket{\bar{j}}{l}$ can now be obtained
by solving an eigenvalue problem in the $\ket{\bar{j}}$-basis:
\begin{align}
	\bra{\bar{j}} L - \lambda_l \ket{l}
		&= \sum_{j'} \frac{\bra{\bar{j}} L - \lambda_l
\ket{\bar{j}'}}{\braket{\bar{j}'}{\bar{j}'}}
				\braket{\bar{j}'}{l}, \\
	0	&= \sum_{j'} \left(\frac{\bar{L}_{jj'}}{\braket{\bar{j}'}{\bar{j}'}}
				- \lambda_l\delta_{jj'}\right) \braket{\bar{j}'}{l}, \\
\label{tet-eig-rec}
	0 &= \sum_{j'} \left(\bar{L}_{jj'} - \lambda_l \braket{\bar{j}'}{\bar{j}} \right)
				\frac{\braket{\bar{j}'}{l}}{\braket{\bar{j}'}{\bar{j}'}},
\end{align}
where $\lambda_l = \lambda(a,b,l)$. Since $\bar{L}_{jj'}$ is
tridiagonal, we obtain a three-term recurrence relation for the
$\braket{\bar{j}}{l}$ transition coefficients. Expanding the expression for
$\bar{L}_{jj'}$, we find the following general form of the recurrence
relation:
\begin{equation}%
\label{tet-rec}
	\frac{\bar{L}_{j,j-2}}{\braket{\overline{j-2}}{\overline{j-2}}}
			\braket{\overline{j-2}}{l}
	+ \left(\frac{\bar{L}_{jj}}{\braket{\bar{j}}{\bar{j}}}-\lambda_l\right)
			\braket{\bar{j}}{l} 
	+ \frac{\bar{L}_{j,j+2}}{\braket{\overline{j+2}}{\overline{j+2}}}
			\braket{\overline{j+2}}{l}
	= 0,
\end{equation}
with the provision that $\bar{L}_{jj'}$ vanishes whenever either of the
indices fall outside the admissible range or $j=j'=0$. Finally, the transition
coefficients are uniquely determined (up to sign) by requiring the
normalization condition
\begin{equation}
	\sum_j \frac{\braket{l}{\bar{j}}\braket{\bar{j}}{l}}
		{\braket{\bar{j}}{\bar{j}}} = \braket ll.
\end{equation}
Practically, it is more convenient to recover the correct normalization
for all $j$ and fixed $l$, or vice versa, by requiring
$\braket{\bar{j}}{l}$ to agree with~\eqref{tet-formula} for $j=\uj$,
cf.~\eqref{jmin}, where the sum reduces to a single term.

Once the $\Tet$-network has been evaluated recursively, the $6j$-symbol
can be obtained from Equation~\eqref{KL-6j}.  Alternatively, a linear,
three-term recurrence relation directly for the $6j$-symbol follows
from~\eqref{tet-rec} and the linear, two-term recurrence relations for
the bubble and $\theta$-networks, obvious from~\eqref{bub-def}
and~\eqref{theta-def}. However, because of the additional normalization
factors in Equation~\eqref{KL-6j}, this direct recurrence relation cannot be
cast in the form of a symmetric eigenvalue problem
like~\eqref{tet-eig-rec} using rational operations alone.

\appendix
\section{Formulas}
For a complex number $q\ne0$ and an integer $n$ the corresponding
\emph{quantum integer} is defined as
\begin{equation}\label{qi-def}
	[n] = \frac{q^n-q^{-n}}{q-q^{-1}} .
\end{equation}
In the limit $q\to 1$, we recover the regular integers, $[n]\to n$.
When $q=\exp(i\pi/r)$, for some integer $r>1$, it is a \emph{primitive root of
unity} and the definition reduces to
\begin{equation}
	[n] = \frac{\sin(n\pi/r)}{\sin(\pi/r)} ,
\end{equation}
This expression is clearly real and positive in the range $0<n<r$.
\emph{Quantum factorials} are direct analogs of classical factorials:
\begin{equation}
	[0]! = 1, \quad [n]! = [1][2]\cdots[n] .
\end{equation}

Next, we give the evaluations of some spin networks needed in the paper.
They are reproduced from Ch.~9 of~\cite{KL}.
The \emph{bubble diagram} evaluates to
\begin{equation}\label{bub-def}
	\loopj = (-)^{j}[j+1]
\end{equation}
whenever it is non-vanishing. For generic $q$, it vanishes if $j<0$ and
if $q$ is a primitive root of unity then it also vanishes when $j>r-2$.
The \emph{$\theta$-network} evaluates to
\begin{equation}\label{theta-def}
	\theta(a,b,c) = \frac{(-)^s [s+1]![s-a]![s-b]![s-c]!}{[a]![b]![c]!} ,
\end{equation}
with $s=(a+b+c)/2$, whenever the twice-spins $(a,b,c)$ are admissible and
vanishes otherwise. Admissibility consists of the following
criteria (besides the obvious $a,b,c\ge0$):
\begin{align}
\label{tri-ineq}
	\text{\emph{triangle inequalities}} ~~ &
	\left\{\begin{aligned}
		a &\le b+c \\
		b &\le c+a \\
		c &\le a+b
	\end{aligned}\right. , \\
\label{parity}
	\text{\emph{parity}} ~~ &
		a+b+c \equiv 0 \pmod{2} .
\intertext{
When $q$ is a primitive root of unity, further constraints needs to be
satisfied:
}
\label{r-bound}
	\text{\emph{$r$-boundedness}} ~~ &
	\left\{\begin{aligned}
		a,b,c &\le r-2 \\
		a+b+c &\le 2r-4
	\end{aligned}\right. .
\end{align}

The \emph{tetrahedral-} or \emph{$\Tet$-network} evaluates to
\begin{equation}\label{tet-formula}
	\Tet(a,b,c,d;j,l)
	= \frac{\mathcal{I}!}{\mathcal{E}!}
		\sum_S
			\frac{(-)^S [S+1]!}{\prod_\imath[S-a_\imath]!\prod_\jmath[b_\jmath-S]!},
\end{equation}
where the summation is over the range $m\le S\le M$ and
\begin{align}
	\mathcal{I}! &= \prod_{\imath,\jmath}[b_\jmath-a_\imath]!, &
		\mathcal{E}! &= [a]![b]![c]![d]![j]![l]!, \\
	a_1 &= (a+d+j)/2, & b_1 &= (b+d+j+l)/2, \\
	a_2 &= (b+c+j)/2, & b_2 &= (a+c+j+l)/2, \\
	a_3 &= (a+b+l)/2, & b_3 &= (a+b+c+d)/2, \\
	a_4 &= (c+d+l)/2, & m &= \max\{a_{\imath}\}, \\
	    &             & M &= \min\{b_{\jmath}\}.
\end{align}
The indices $\imath$ and $\jmath$ fully span the defined ranges.
Each of the triples of twice-spins $(a,b,l)$, $(c,d,l)$, $(a,d,j)$ and
$(c,b,j)$ must be admissible, otherwise the $\Tet$-network vanishes.
Then, due to parity constraints, $a_{\imath}$, $b_{\jmath}$,
$m$, $M$, and $S$ are always integers. If the twice-spins $(a,b,c,d)$
are fixed, the admissibility conditions for generic $q$ enforce
the ranges of $\uj\le j\le\oj$ and $\ul\le l\le\ol$ to
\begin{align}
\label{jmin}
	\uj &= \max\{|a-d|,|b-c|\}, & \oj &= \min\{a+d,b+c\}, \\
\label{lmin}
	\ul &= \max\{|a-b|,|c-d|\}, & \ol &= \min\{a+b,c+d\},
\end{align}
with
\begin{align}
	j &\equiv a+b \equiv c+d \pmod{2}, \\
	l &\equiv a+d \equiv b+c \pmod{2}.
\end{align}
The number of admissible values is the same for $j$ and $l$ and is equal
to $n=\max\{0,\bar{n}\}$, where
\begin{align}
\label{n-dim}
	\bar{n}&=\min\{m,s-M\}+1,  & m&=\min\{a,b,c,d\}, \\
	s      &=(a+b+c+d)/2,      & M&=\max\{a,b,c,d\}.
\end{align}
This number $n$ is also the dimension of the space of $4$-valent spin
networks with fixed twice-spins $(a,b,c,d)$ labeling the free edges.
This dimension is non-vanishing, $n>0$, precisely when the twice-spins
satisfy the conditions
\begin{align}
\label{abcdmax}
	a+b+c+d &\le 2\max\{a,b,c,d\}, \\
\label{abcdpar}
	a+b+c+d &\equiv 0 \pmod{2}.
\end{align}

When $q$ is a primitive root of unity, the admissible ranges shrink to
$\uj\le j\le\oj_r$ and $\ul\le l\le\ol_r$, where
\begin{align}
	\oj_r &= \min\left\{\oj,r-2,2r-4-\max\{a+d,b+c\}\right\}, \\
	\ol_r &= \min\left\{\ol,r-2,2r-4-\max\{a+b,c+d\}\right\}.
\end{align}
The number of admissible values in each range is thus restricted to
$n=\max\{0,\bar{n}_r\}$, where
\begin{equation}
\label{nr-dim}
	\bar{n}_r = \min\left\{\bar{n}, r-1-\max\{M,s-m\}\right\}.
\end{equation}
The condition $n>0$ requires~\eqref{abcdmax}, \eqref{abcdpar} and
\begin{equation}
\label{abcdrbd}
	a+b+c+d \le 2\min\{a,b,c,d\}+2r-4 .
\end{equation}

The above admissibility criteria are well known. However, the consequent
explicit expressions for the constraints on $(a,b,c,d)$, the bounds on
$j$ and $l$, and the dimension $n$ are not easily found in the
literature.

In the classical $q=1$ case, the Kauffman-Lins version of the
$6j$-symbol~\eqref{KL-6j} differs from the Racah-Wigner convention used
in the physics literature, which preserves the symmetries of the
underlying $\Tet$-network. The two $6j$-symbols are related through the
formula
\begin{equation}%
\label{RW-KL}
	\sixjrw{j_1/2}{j_2/2}{j_3/2}{J_1/2}{J_2/2}{J_3/2} = 
	\frac{\Tet(J_1,J_2,j_1,j_2;J_3,j_3)}
		{\sqrt{|\theta(J_1,J_2,j_3)\theta(j_1,j_2,j_3)
			\theta(J_1,j_2,J_3)\theta(J_2,j_1,J_3)}|},
\end{equation}
which can be obtained by comparing the explicit
expressions~\eqref{tet-formula} and those of Sec.~9.2.1 of~\cite{VMKh}.
Note that the argument of the absolute value under the square root has
sign $(-)^{j_3-J_3}$. 

\section*{Acknowledgements}
The author thanks Dan Christensen for many helpful discussions as well
as Nicolai Reshetikhin and an anonymous referee for pointing out
previous work on the above recurrence relation. In the course of this
work, the author was supported by Postgraduate (PGS) and Postdoctoral
(PDF) Fellowships from the Natural Science and Engineering Research
Council (NSERC) of Canada.

\newpage
\bibliographystyle{ws-ijgmmp}
\bibliography{eigenrec}

\end{document}